\documentclass{article}

\usepackage{amssymb}
\usepackage{amsfonts}
\usepackage{latexsym}
\usepackage[T1,T2A]{fontenc}
\usepackage[cp1251]{inputenc}
\usepackage[T2B]{fontenc}
\usepackage[bulgarian,russian,english]{babel}

\newtheorem {theor} {\bf Theorem}

\newtheorem {zad} {\bf Problem}
\newtheorem{crl} {\bf Corollary}

\title {ON A CLASS OF BINARY MATRICES}
\author {Krasimir Yankov Yordzhev}
\date {}

\begin {document}
\inputencoding{cp1251}

\maketitle
\begin{abstract} The paper studies the set of all square binary
matrices containing an exact number of 1's in each rows and in
each column. A connection is established between the cardinal
number of this set and the cardinal number of its  subset of
matrices containing 1 in the lower right corner. With the help of
this result a new proof is  advanced of the I. Good and J. Grook
theorem. In connection with the firs result a classification has
also been made of square binary matrices containing three 1's in
each row and column and 1 in the lower right corner.
\end{abstract}

\section{Introduction} A {\it binary} (or  {\it boolean}, or {\it (0,1)-matrix})
 is a matrix whose all elements belong to the set ${\bf B}
=\{ 0,1 \}$. With ${\bf B}_n$ we will denote the set of all  $n
\times n$  binary matrices.

Using the notation from   \cite {1}, we will call
$\Lambda_n^k$-matrices all $n\times n$ binary matrices in each row
and each column of which there are exactly $k$ in number 1's.

Let us accept the associative operations $"+"$ and $"\cdot " $  in
the set ${\bf B} =\{ 0,1 \}$ defined as: $0+0=0,$ $1+0=0+1=1+1=1,$
$0\cdot 0=1\cdot 0=0\cdot 1=0,$ $1\cdot 1=1$. Let the scalar
product of two vectors (defined by the n-tuple) with elements from
${\bf B}$ be overdetermined by the operations thus introduced. Let
us accept that in ${\bf B_n }$ the common definition for the
operation of matrix product as a scalar product of the
corresponding row vector and column vector is over determined by
the operation of scalar product introduced above. Considering all
this,  ${\bf B}_n$, along with the operation matrix product, is a
semigroup that is isomorphic with the semigroup ${\cal B}(M)$
comprising all binary relation in a  set  $M$, where
$|M|=n<\infty$
  (see for example \cite{lallement}).

Let us consider the following combinatorial problem:

\begin{zad} \label{z1}
Find out the number of all binary relations  $\omega \in {\cal
B}(M)$, $|M|=n<\infty$, such that for all $a\in M$ the equation
$$|\{ x\in M\; |\; (a,x)\in \omega \} |=|\{ x\in M\; |\;  (x,a) \in \omega \} |=k ,$$
is correct, where  $n$ and $k$ are natural numbers.
\end{zad}

Using the language of graph theory  (see for example \cite{aigner}
or \cite{harary}) problem \ref{z1} is equivalent to

\begin{zad}\label{z2}  Find out the number of all bipartite graphs $G=(V_1 \cup V_2 ,E)$,
such that $|V_1 |=|V_2 |=n$ and each vertex is incident with
exactly $k$ edges.
\end{zad}

Clearly, problems  \ref{z1} and \ref{z2} can be reduced  to
solving the following

\begin{zad}\label{z3} Find out the number of all   $n\times n$
matrices containing exactly $k$ 1's in each row and each column,
e.g. the number of all $\Lambda_n^k$-matrices.
\end{zad}

The goal of this paper is to consider certain special cases of
problem  \ref{z3}.

Let us denote the number of all $\Lambda_n^k$-matrices with
$\lambda_{n,k}$.

In \cite{6} is offered the formula:

\begin{equation}\label{1}
\lambda_{n,2} = \sum_{2x_2 +3x_3 + \cdots +nx_n =n}
\frac{(n!)^2}{\displaystyle \prod_{r=2}^n x_r !(2r)^{x_r}}
\end{equation}

One of the first recursive formulas for the calculation of
$\lambda_{n,2}$ appeared in \cite{anand} :
\begin{equation}\label{2}
\begin{array}{l}
\displaystyle \lambda_{n,2} =\frac{1}{2} n(n-1)^2 \left[ (2n-3)
\lambda_{(n-2),2} +(n-2)^2 \lambda_{(n-3),2} \right] \quad \mbox{for} \quad n\ge 4 \\
\displaystyle \lambda_{1,2} =0,\quad \lambda_{2,2} =1,\quad
\lambda_{3,2} =6
\end{array}
\end{equation}

Another recursive formula for the calculation of $\lambda_{n,2}$
occurs  in \cite{13} :
\begin{equation}\label{3}
\begin{array}{l}
\displaystyle \lambda_{n,2} = (n-1)n\lambda_{(n-1),2}
+\frac{(n-1)^2
n}{2} \lambda_{(n-2),2} \quad \mbox{for} \quad n\ge 3 \\
\lambda_{1,2} =0,\quad \lambda_{2,2} =1
\end{array}
\end{equation}

The following recursive system  for the calculation of
$\lambda_{n,2}$ is put forward in \cite{iord}:
\begin{equation}\label{4}
\left|
\begin{array}{l}
\lambda_{(n+1),2} =n(2n-1)\lambda_{n,2} + n^2 \lambda_{(n-1),2} - \pi_{n+1} \quad \mbox{for} \quad n\ge 2\\
 \pi_{n+1} = {n^2 (n-1)^2 \over 4} [8(n-2)(n-3) \lambda_{(n-2),2} +(n-2)^2 \lambda_{(n-3),2} -4\pi_{n-1}]\quad \mbox{for} \quad n\ge 4 \\
\lambda_{1,2} = 0 ,\quad \lambda_{2,2} = 1 ,\quad \pi_1 = \pi_2 =
\pi_3 = 0 ,\quad \pi_4 = 9
\end{array}
\right.
\end{equation}
where  $\pi_n$ identifies the number of a special class of
$\Lambda_n^2$-matrices.

For the classification  of all non defined concepts and notations
as well as for common assertion which have not been proved here,
we recommend sources \cite{aigner,lankaster,baranov,tar2}.

\section{On a partition of the set  $\Lambda_n^k$}\label{sect1}
Let us introduce the notations:
\begin{equation}\label{L+}
\Lambda_{n}^{k+} =\vert \{A=(a_{ij} )\in \Lambda_{n}^k \vert \;
a_{n\, n} =1 \}\vert
\end{equation}
\begin{equation}\label{L-}
\Lambda_{n}^{k-} =\vert \{A=(a_{ij} )\in \Lambda_{n}^k \vert \;
a_{n\, n} =0 \}\vert
\end{equation}

Obviously
\begin{equation}\label{part} \Lambda_{n}^{k+} \cap
\Lambda_{n}^{k-} =\emptyset \quad \mbox{and} \quad
\Lambda_{n}^{k+} \cup \Lambda_{n}^{k-} =\Lambda_{n}^{k},
\end{equation}
in other words $\{ \Lambda_n^{k+} ,\Lambda_n^{k-} \}$ represents a
partition of the set  $\Lambda_n^k$.

We set:

\begin{equation}\label{lpl}
\lambda_{n,k}^{+} =|\Lambda_{n}^{k+} |
\end{equation}
\begin{equation}\label{lmn}
\lambda_{n,k}^{-} =|\Lambda_{n}^{k-} |
\end{equation}

Formula   (\ref{3}) occurs for the first time in  \cite{13} and it
has been deduced in a manner applicable only to the calculation of
the number of the $\Lambda_n^2$-matrices. The method for the
obtaining of the recursive relation (\ref{3}) which  we offer and
which we will describe in section \ref{gg} will being as closer to
the discovery  of the analogical formula for values greater than
$k$. In this case  $k$ represents the number of units in each row
and each column of the respective square matrices. The method is
grounded in the following assertion:

\begin{theor} \label{lema1} It is true the equation

\begin{equation}\label{lll}
\lambda_{n}^{k-} =\frac{n-k}{k} \lambda_{n}^{k+} ,
\end{equation}
where  $\lambda_{n,k}^{+} $ and $ \lambda_{n,k}^{-}$ are set by
formulas     (\ref{lpl}) and (\ref{lmn}) respectively.
\end{theor}

Proof. Let us accept that $A$ and $B$ are
$\Lambda_{n}^k$-matrices. We will say that $A$ and $B$ are
$\rho$-equivalent   $( A\rho B )$, if the removing of the columns
ending in 1 results in equal $n\times (n-k)$ matrices. Obviously,
$\rho$ is an equivalence relation. We use $\rho_A$ to denote the
set of elements to which $A$ is related by the equivalence
relation $\rho$.

Let $A=(a_{ij} )$ be a $\Lambda_{n}^k$-matrix. Let's denote with
$p^+$ the number of all matrices  $\rho$-equivalent to $A$ in
which the element in the lower right corner is equal to 1 and with
 $p^-$ the number of all matrices $\rho$-equivalent to   $A$ in which
 the element in the lower right
 corner is equal to  0. Let $K_{j_1} ,K_{j_2} ,..., K_{j_k}$ are
the row-vectors of matrix  $A$ with 1 in final position. The set
$J=\{ j_1 ,j_2 ,..., j_k \}$ is partitioned into subsets $J_{r}$,
$r=1,2,\ldots s$, such that $j_u$ and $j_v$ are part of the same
subset if and only if $K_{j_u} =K_{j_v}$. It is easy to detect
that $\displaystyle J=\bigcup_{r=1}^s J_{r}$ and $J_{u} \cap J_{v}
=\emptyset$ for $u\ne v$. We set $k_r =|J_r |$, $r=1,2,\ldots ,s$.
Obviously
\begin{equation}\label{kkk}
k_1 +k_2 +\cdots +k_s =k.
\end{equation}

Let  $C$ be the $n\times (n-k)$ matrix which comes from  $A$ by
removing the columns  $K_{j_1} ,K_{j_2} ,\ldots , K_{j_k} $. In
this case, with the help of the different ways of adding new
columns to $C$, we will obtain all elements of the set $\rho_A$.
Let us first we add $k_1$ columns which  equal to those columns of
$A$ whose numbers belong to the set $J_1$. This can be done using
${n-k+k_1 \choose k_1}$ number of ways. We can then add $k_2$
equal columns in ${n-k+k_1 +k_2 \choose k_2}$ possible ways. These
equal columns are also equal to the columns in $A$ with number
tags from  $J_2$, etc. Therefore

\medskip
$\displaystyle \vert \rho_A \vert ={{n-k+k_1}\choose k_1}{{n-k+k_1
+k_2 }\choose k_2 } \cdots {{n-k+k_1 +k_2 +\cdots +k_s}\choose
k_s} =$
\medskip

\noindent $\displaystyle  ={{n}\choose k_s} {{n-k_s }\choose
k_{s-1}} {{n-k_s -k_{s-1}}\choose k_{s-2}} \cdots {{n-k_s -k_{s-1}
-k_{s-2} +\cdots +k_2}\choose k_1} =$
\medskip

\noindent$\displaystyle =\frac{n!(n-k_s )!(n-k_s -k_{s-1} )!\cdots
(n-k_s -k_{s-1} -\cdots -k_2 )!}{k_s !(n-k_s )! k_{s-1} ! (n -k_s
-k_{s-1} )!  \cdots k_1 ! (n -k_s -k_{s-1} -\cdots -k_2 -k_1 )!}
=$

\noindent$\displaystyle = \frac{n!}{k_1 !k_2 !\cdots k_s !(n-k)!}$

\medskip
Analogically, for $p^-$ we get

\medskip
$\displaystyle p^- =\frac{(n-1)!}{k_1 !k_2 ! \cdots k_s !
(n-1-k)!}$,\\
having the mind the fact that  we cannot add new columns after the
last column of matrix $C$.

For  $p^+$ we obtain the equation:

 $\displaystyle p^+ = \vert \rho_A \vert -p^- =
\frac{n!}{k_1 !k_2 ! \cdots k_s ! (n-k)!} -\frac{(n-1)!}{k_1 !k_2
! \cdots k_s ! (n-1-k)!} =$

\noindent$\displaystyle = \frac{k(n-1)!}{k_1 !k_2 ! \cdots k_s !
(n-k)!}$

\medskip Then $\displaystyle \frac{p^-}{p^+} =\frac{n-k}{k}$,
e.g. $\displaystyle p^- =\frac{n-k}{k} p^+$. Summarize by
equivalence classes we arrive a the equation we were supposed to
prove.

Considering equations  (\ref{L+}) $\div$ (\ref{lmn}) and theorem
\ref{lema1}, it is possible to formulate
\begin{crl}\label{col}\it
\begin{equation}\label{sl1}
\lambda_{n,k} =\lambda_{n,k}^+ +\lambda_{n,k}^- =\lambda_{n,k}^+
+\frac{n-k}{k}\lambda_{n,k}^+ =\frac{n}{k}\lambda_{n,k}^+
\end{equation}
\end{crl}

\section{ Some applications}

Theorem \ref{lema1} and corollary  \ref{col} are useful in that
they facilitate the calculation $\Lambda_n^{k+}$ as compared to
that of all  $\Lambda_n^k$-matrices.

\subsection{ A different proof of the  {\bf  I. Good} and {\bf J. Grook } theorem}\label{gg}

Using corollary \ref{col} in order to obtain a formula for the
 $n\times n$ binary matrices, it is enough to find a formula for $\lambda_{n,2}^+$. We can
this with the help of the following

\begin{theor}\label{t2} If  $n\ge 3$, then
$$\lambda_{n,2}^+  = 2(n-1)\lambda_{(n-1),2} +(n-1)^2 \lambda_{(n-2),2} $$
\end{theor}

\par Proof. Let $A=(a_{ij})$ be a $\Lambda_{n-1}^2$-matrix. The
matrix  $A$ can be give rise to the  $\Lambda_{n}^2$-matrix
$B=(b_{ij})$ following manner: We choose $p,q$ such that   $a_{p\,
q} =1$. This can be accomplished in $2(n-1)$ ways. We set $b_{p\,
q} =0 $, $b_{p\, n} =b_{n\, q} =b_{n\, n} =1 $, $b_{i\, n} =b_{n\,
j} =0$  and $b_{ij} =a_{ij}$ for $1\le i,j\le n-1$, $i\ne p$,
$j\ne q$. It is easy to see that $B$ s $\Lambda_{n}^2$-matrix with
1 in the lower right corner. Besides $p$ and $q$ can be identified
uniquely through  $B$ and matrix $A$ can be restored. Consequently
$\lambda_{n,2}^+ =2(n-1)\lambda_{(n-1),2} +t$, where $t$ is the
number of all $\Lambda_{n}^2$-matrices containing 1 in the lower
right corner, which cannot be generated in the manner described
above. These are $\Lambda_{n}^2$-matrices,  $B=(b_{ij})$ for which
there are  $p$ and $q$ such that   $b_{p\, q} =b_{n\, q} =b_{p\,
n} =b_{n\, n} =1$ and these are the only 1's (2 in each row and
column) in rows with number $p$ and  $n$ and in columns with
number $q$ and $n$. In this case, however, removing rows number
$p$ and $n$ and columns number $q$ and  $n$, we obviously obtain a
$\Lambda_{n-2}^2$-matrix. On the contrary, each
$\Lambda_{n-2}^2$-matrix can give rise to a $\Lambda_{n}^2$-matrix
by inserting two  new rows, their numbers will be $p$ and $n$ and
two new columns, their numbers will be $q$ $n$, with 0 everywhere
except for the places  of intersection. Since  $p$ and $q$ vary
from $1$ to $n-1$, then    $t=(n-1)^2 \lambda_{(n-2),2}$. This
proves the theorem.

Applying theorems \ref{lema1} and  \ref{t2} directly we obtain:

\begin{theor}\cite{13}\label{good} The number of all $n\times n$ square
binary matrices with exactly two 1's in each row and each column
is given by the next formula:
$$\begin{array}{l}
\displaystyle \lambda_{n,2} = (n-1)n\lambda_{(n-1),2}
+\frac{(n-1)^2
n}{2} \lambda_{(n-2),2} \quad \mbox{for} \quad n\ge 3\\
\lambda_{1,2} =0,\quad \lambda_{2,2} =1
\end{array} $$
\end{theor}

\subsection{ On the number of $\Lambda_n^3$-matrices}

The following formula in an explicit form for the calculation of
$\lambda_3 (n)$ is offered in \cite{13}.

\begin{equation}\label{111}
 \lambda_3 (n)=\frac{n!^2}{6^n} \sum \frac{(-1)^\beta (\beta
+3\gamma )! 2^\alpha 3^\beta}{\alpha !\beta ! \gamma !^2 6^\gamma}
\end{equation}
where the sum is done as regard all $ \frac{(n+2)(n+1)}{2}$
solutions in nonnegative whole numbers of the equation   $\alpha
+\beta+\gamma =n$.

As it is noted in  \cite{stan} formula    (\ref{111}) does not
give us good opportunities to study behavior of  $\lambda_{n,3}$.
The aim of the current consideration is to go one step closer to
the obtaining of a new recursive formula for the calculation of
$\lambda_{n,3}$, which could help avoid certain inconveniences
resulting from use of formula (\ref{111}).

Let  $X=(x_{ij} )\in \Lambda_{n}^{3+} $ and all 1's in the last
columns and in the last row are respectively  the elements $x_{s\,
n} ,$ $x_{t\, n} ,$ $x_{n\, p} ,$ $x_{n\, q} ,$ $x_{n\; n}$, where
$s,t,p,q\in \{ 1,2,\ldots ,n-1\} ,$ $s\ne t,$ $p\ne q$.
$\widetilde{X}$  will be denoted $  2\times 2$ submatrix

\begin{equation}\label{Xtilde}
\widetilde{X} =\left(
\begin{array}{cc}
x_{s\, p} & x_{s\, q} \\
x_{t\, p} & x_{t\, q}
\end{array}
\right)
\end{equation}

The set  $\Lambda_{n}^{3+}$ is partitioned into the following
nonintersecting subsets:
\begin{eqnarray}
{\rm A}_{n} &=& \left\{ X\in \Lambda_{n}^{3+} \left| \widetilde{X}
=\left(
\begin{array}{cc}
1 & 1\\
1 & 1\\
\end{array}
 \right) \right. \right\} \\
{\rm B}_{n} &=&\left\{ X\in \Lambda_{n}^{3+} \left| \widetilde{X}
\in\left\{ \left(
\begin{array}{cc}
0 & 1\\
1 & 1\\
\end{array}\right) ,\left(
\begin{array}{cc}
1 & 0\\
1 & 1\\
\end{array}\right) , \left(
\begin{array}{cc}
1 & 1\\
0 & 1\\
\end{array}\right) , \left(
\begin{array}{cc}
1 & 1\\
1 & 0\\
\end{array}\right) \right\} \right. \right\} \\
\Gamma_{n} &=&\left\{ X\in \Lambda_n^{3+} \left| \widetilde{X}
\in\left\{ \left(
\begin{array}{cc}
1 & 0\\
1 & 0\\
\end{array} \right) ,\left(
\begin{array}{cc}
0 & 1\\
0 & 1\\
\end{array}\right)\right\}  \right. \right\} \\
\Delta_{n} &=&\left\{ X\in \Lambda_n^{3+} \left| \widetilde{X}
\in\left\{ \left(
\begin{array}{cc}
1 & 1\\
0 & 0\\
\end{array} \right) ,\left(
\begin{array}{cc}
0 & 0\\
1 & 1\\
\end{array}\right)\right\}  \right. \right\} \\
{\rm E}_{n} &=&\left\{ X\in \Lambda_n^{3+} \left| \widetilde{X}
\in\left\{ \left(
\begin{array}{cc}
1 & 0\\
0 & 1\\
\end{array} \right) ,\left(
\begin{array}{cc}
0 & 1\\
1 & 0\\
\end{array}\right)\right\}  \right. \right\} \\
{\rm Z}_{n} &=&\left\{ X\in \Lambda_n^{3+} \left| \widetilde{X}
\in\left\{ \left(
\begin{array}{cc}
1 & 0\\
0 & 0\\
\end{array} \right) ,\left(
\begin{array}{cc}
0 & 1\\
0 & 0\\
\end{array}\right) ,
\left(
\begin{array}{cc}
0 & 0\\
1 & 0\\
\end{array} \right) ,\left(
\begin{array}{cc}
0 & 0\\
0 & 1\\
\end{array}\right)\right\}  \right. \right\} \\
{\rm H}_{n} &=&\left\{ X\in \Lambda_n^{3+} \left| \widetilde{X}
=\left(
\begin{array}{cc}
0 & 0\\
0 & 0\\
\end{array} \right)  \right. \right\}
\end{eqnarray}

We set
\begin{equation}\label{alphan}
\begin{array}{c}
\alpha_n =| {\rm A}_{n} |,\quad \beta_n=|{\rm B}_{n} |,\quad
\gamma_n=| \Gamma_{n} |,\quad  \delta_n=|\Delta_{n} |,\\
\epsilon_n =| {\rm E}_{n} |,\quad \zeta_n=| {\rm Z}_{n} |,\quad
\eta_n=|{\rm H}_{n} |
\end{array}
\end{equation}

Obviously
\begin{equation}\label{gamman}
\gamma_n=\delta_n
\end{equation}
and
\begin{equation}\label{suma}
\lambda_{n,3}^+ =\alpha_n +\beta_n +\gamma_n +\delta_n +\epsilon_n
+\zeta_n +\eta_n .
\end{equation}

\begin{theor}\label{t4}
$$\lambda_{n,3}^+ =\frac{3(n-1)(3n-8)}{2} \lambda_{(n-1),3} +\alpha_n +\beta_n +2\gamma_n -\eta_n $$
\end{theor}

Proof. Let  $Y=(y_{i\, j} )\in\Lambda_{n-1}^{3}$. In $Y$ we choose
two 1's not belonging to the same row or column. Let then these be
the elements  $y_{s\, p}$ and $y_{t\, q} $, $s,t,p,q\in  \{
1,2,\ldots ,n-1\} ,$ $s\ne t,$ $p\ne q$. This can happen in
$\frac{3(n-1)[3(n-1)-5]}{2} = \frac{3(n-1)(3n-8)}{2}$ ways. The
1's thus selected are turned into 0's and in  $Y$ in the last
place  one more column  (number $n$) and one more row (number $n$)
are added, so that  $y_{s\, n} =y_{t\, n} =y_{n\, p} =y_{n\, q}
=y_{n\, n} =1$ и $y_{i\, n} =y_{n\; j} =0$ for $i\notin \{
s,t,n\},$ $j\notin \{ p,q,n\}$. Obviously, the matrix thus formed
belongs to one of the set  ${\rm E}_{n} ,$ ${\rm Z}_{n}$ or ${\rm
H}_{n}$.

On the contrary, let $X=(x_{ij} )$ be a matrix from ${\rm E}_{n}$
or ${\rm Z}_{n}$.  $\widetilde{X}$ then has unique zero diagonal
whose elements we turn into 1's and remove the last row (number
$n$) and the last column (number $n$). In this way a
$\Lambda_{n-1}^3$-matrix is generated.

Let  $X =(x_{ij} )\in {\rm H}_{n} $ and  $\widetilde{X} = \left(
\begin{array}{cc}
x_{sp} & x_{sq}\\
x_{tp} & x_{tq}
\end{array} \right) =
\left(
\begin{array}{cc}
0 & 0\\
0 & 0\\
\end{array} \right)$.
We select one   diagonal of $\widetilde{X}$ and turn 0's of this
diagonal  into 1's. Remove the last row $n$ and the last column
$n$. We get a $\Lambda_{n-1}^3$-matrix. Consequently for each  ${\rm
H}_{n}$-matrix correspond two  $\Lambda_{n-1}^3$-matrices. Then,
$$\lambda_{n,3}^+ =\frac{3(n-1)(3n-8)}{2} \lambda_{(n-1),3} -\eta_n+t$$
where $t=\alpha_n+\beta_n+\gamma_n+\delta_n$.  Considering
 (\ref{gamman}) and (\ref{suma}) we proved the theorem.

Krasimir Yankov Yordzhev

South-West University  ''N. Rilsky''

2700 Blagoevgrad

Bulgaria

\rm e-mail:  iordjev@yahoo.com, iordjev@aix.swu.bg

\begin {thebibliography}{99}
\bibitem{aigner} \textsc{M. Aigner} Combinatorial theory. Springer-Verlag, 1979.
\bibitem{anand}  \textsc{H. Anand,  V. C. Dumir, H. Gupta} A combinatorial distribution problem.
\textit{Duke Math. J.} 33 (1966), 757-769.
\bibitem{13}  \textsc{I. Good, J. Grook} The enumeration of arrays and generalization related to contingency tables.
\textit{Discrete Math}, 19 (1977), 23-45.
\bibitem{harary} \textsc{F. Harary} Graph theory. Addison-Wesley,
1969.
\bibitem{lallement} \textsc{G. Lallement} Semigroups and combinatorial application. J. Wiley \& Sons, 1979.
\bibitem{lankaster} \textsc{P. Lancaster} Theory of Matrices. Academic Press, NY, 1969.
\bibitem{stan} \textsc{R. P. Stanley} Enumerative combinatorics. V.1, Wadword \& Brooks, California, 1986.
\bibitem{baranov} \textsc{В. И. Баранов, Б. С. Стечкин} Экстремальные номбинаторные задачи и их приложения.
Москва, Наука, 1989
\bibitem{iord} \textsc{К. Я. Йорджев} Комбинаторни задачи над бинарни матрици.
\textit{Математика и математическо образование},  24 (1995),
288-296.
\bibitem{6} \textsc{В. Е. Тараканов} Комбинаторные задачи на бинарных матрицах. \textit{Комбинаторный анализ},
Москва, изд-во МГУ, 1980, вып.5, 4-15.
\bibitem{tar2} \textsc{В. Е. Тараканов} Комбинаторные задачи и (0,1)-матрицы. Москва, Наука, 1985.
\bibitem{1} \textsc{В. С. Шевелев} Редуцированные латинские прямоугольники и квадратные матрицы с одинаковыми
суммами в стоках и столбцах. \textit{Дискретная математика}, том
4, вып. 1, 1992, 91-110.
\end{thebibliography}
\end{document}